\documentclass[12pt]{article}
\usepackage{amsthm,amsfonts, amsbsy, amssymb,amsmath,graphicx}
\usepackage{graphics}

\author{D.Yu.Krylov, V.O.Manturov}

\date{}

\title{Parity and Relative Parity in Knot Theory}

\begin{document}

\maketitle

\begin{abstract}

In the present paper we give a simple proof of the fact that the set
of virtual links with orientable atoms is closed. More precisely,
the theorem states that if two virtual diagrams $K$ and $K'$ have
orientable atoms and they are equivalent by Reidemeister moves, then
there is a sequence of diagrams $K = K_1 \to \dots \to K_n=K'$ all
having orientable atoms where $K_i$ is obtained from $K_{i-1}$ by a
Reidemeister move. The initial proof heavily relies on the topology
of virtual links and was published in \cite{IM}. Our proof is based
on the notion of parity which was introduced by the second named
author in 2009.

We split the set of crossings of a virtual link diagram into sets of
{\it odd} and {\it even} in accordance with a fixed rule. The rule
must only satisfy several conditions of Reidemeister's type. Then
one can construct functorial mappings of link diagrams by using
parity. The concept of parity allows one to introduce new invariants
and strengthen well-known ones \cite{Ma1}.

\end{abstract}


\section{Parity and Relative Parity}
\

In the present section we introduce the concept of parity. Section 2
contains a short description of the theory of atoms. The main
theorem stated in the Abstract is proved in the Section 3.

{\bf Definition 1.1}~A 4-valent graph $\Gamma$ is called {\it
framed} if for each vertex of $\Gamma$ the four emanating half-edges
are split into two pairs of (formally) opposite.

One can define Reidemeister moves for framed 4-graphs. The first
move is an addition/removal of a loop, see Fig. 1, left. The second
Reidemeister move adds/removes a bigon formed by a pair of edges
which are adjacent in two vertices, see Fig. 1, center. The third
Reidemeister move is shown in Fig. 1, right.


\begin{figure}
\centering\includegraphics[width=200pt]{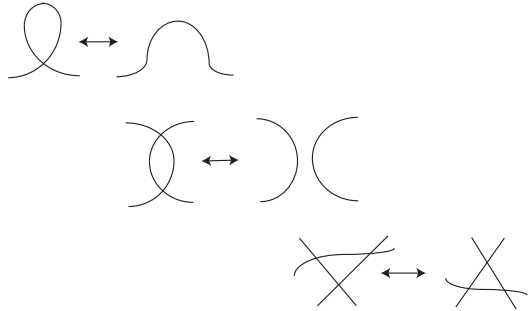} \caption{Three
Reidemeister moves}
\end{figure}

{\bf Definition 1.2}~An equivalence class of framed 4-graphs modulo
Reidemeister moves is called {\it a free link}.

{\bf Definition 1.3}~Two edges $a, a'$ belong to the same {\it
unicursal component} of a graph $\Gamma$ if there exist sequences of
edges $a=a_1, \ldots, a_{n+1}=a'$ and vertices $c_1, \ldots, c_n$
such that for $i=1,\ldots,n$ the edges $a_{i},a_{i+1}$ are opposite
in $c_i$.

The number of unicursal components is an invariant under
Reidemeister moves.

{\bf Definition 1.4}~By {\it a parity} we mean a rule for endowing
all vertices of all 4-valent framed graphs (or all vertices of all
graphs from a subset closed under Reidemeister moves) by elements of
$\mathbb{Z}_2$ such that the following axioms hold:

1.~Each crossing taking part in the first Reidemeister move is
endowed with 0.

2.~Each two crossings taking part in the second Reidemeister move
are either both odd (i.e. are endowed with 1) or both even (i.e. are
endowed with 0).

3.~The sum (modulo 2) of parities of three vertices participating in
the third Reidemeister move equals 0. Moreover, parities of
corresponding vertices ($a$ and $a'$, $b$ and $b'$, $c$ and $c'$ in
Fig. 2) are the same.

Also we require that vertices which do not take part in a move
preserve their parities.

\begin{figure}
\centering\includegraphics[width=160pt]{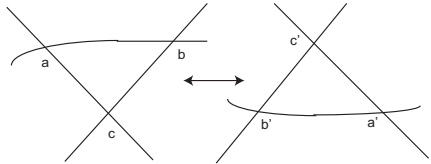}
\caption{Corresponding vertices}
\end{figure}

{\bf Example 1.1}~Free knots (i.e., one-component free links) are
encoded by Gauss diagrams without signs and arrows in such a way
that there is a natural correspondence between vertices of a knot
and chords of its diagram. Two chords are called {\it linked} if
their ends are alternated on the circle of the Gauss diagram. We say
that a chord of a Gauss diagram is even if the number of chords it
is linked with, is {\it even}, and {\it odd} otherwise. It can be
easily checked that this yields a parity. We call it {\it Gaussian
parity}.

{\bf Example 1.2}~This example deals with graphs having 2 unicursal
components. A vertex is even if all four incident half-edges belong
to the same component; otherwise a vertex is odd.

Therefore a parity for the given graph $\Gamma$ is described by a
pair $(\Gamma, \mathcal{P})$. Here $\mathcal{P}$ is a function on
the set of vertices of $\Gamma$ to $\mathbb{Z}_2$ satisfying axioms
1)-3). In the case of more complicated knot-like objects by a {\it
parity} we mean also a rule satisfing the same axioms under
Reidemeister moves.

{\bf Definition 1.5}~{\it A virtual knot diagram} is an image of a
smooth immersion $f\colon S^1\to \mathbb{R}^2$ such that each
intersection point is double, transverse and endowed with classical
(with a choice for underpass and overpass specified) or virtual
crossing structure. {\it A virtual knot} is an equivalence class of
virtual knot diagrams modulo generalized Reidemeister moves. The
notion of a virtual knot was introduced by L.~H.~Kauffman in
\cite{K}.

In the case of virtual link diagrams (where chords of the Gauss
diagram are endowed with signs and arrows) the parity axioms must
hold only for those vertices which can participate in a Reidemeister
move. Note that each parity for free knots induces a parity for
virtual knots, but not vice versa.

Let $\Gamma$ be a framed 4-graph with one unicursal component. The
group $H_1(\Gamma,\mathbb{Z}_2)$ is generated by "halves" which
correspond to vertices: for every vertex $v$ there exist two halves
$\Gamma_{1,v}$ and $\Gamma_{2,v}$ (see Fig. 3) obtained by smoothing
of the graph $\Gamma$ at the vertex $v$. If the set of framed
4-graphs is endowed with a parity, one can consider a cohomology
class $h$ defined by equalities $h( \Gamma_{1,v})=h(
\Gamma_{2,v})=\mathcal{P}(v)$.

\begin{figure}
\centering\includegraphics[width=200pt]{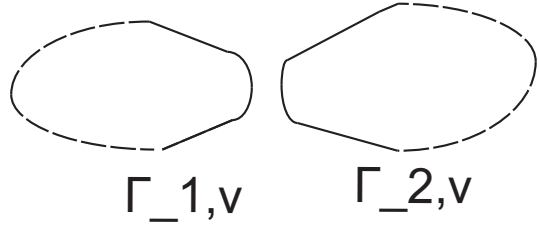} \caption{Halves
$\Gamma_{1,v}$ and $\Gamma_{2,v}$}
\end{figure}

Collecting the properties of this cohomology class we see that:
\begin{enumerate}

\item
 For every framed 4-graph $\Gamma$ we have $h(\Gamma) = 0$.

\item
 If $\Gamma'$ is obtained from $\Gamma$ by a first Reidemeister
move adding a loop then for every basis $\{\alpha_i\}$ of
$H_1(\Gamma, \mathbb{Z}_2)$ there exists a basis of the group
$H_1(\Gamma', \mathbb{Z}_2)$ consisting of one element $\beta$
corresponding to the loop and a set of elements $\alpha_i'$
naturally corresponding to $\alpha_i$. Then we have $h(\beta) = 0$
and $h(\alpha_i) = h(\alpha_i')$.

\item
 Let $\Gamma'$ be obtained from $\Gamma$ by a third increasing
Reidemeister move. Then for every basis $\{\alpha_i\}$ of
$H_1(\Gamma, \mathbb{Z}_2)$ there exists a basis in $H_1(\Gamma',
\mathbb{Z}_2)$ consisting of one "bigon" $\gamma$, the elements
$\alpha_i'$ naturally corresponding to $\alpha_i$, and one
additional element $\delta$, see Fig. 4, left. Then the following
holds: $h(\alpha_i) = h(\alpha_i'), h(\gamma) = 0$.

\item
 Let $\Gamma'$ be obtained from $\Gamma$ by a third Reidemeister move.
Then there exists a graph $\Gamma''$ with one vertex of valency 6
and the other vertices of valency 4 which is obtained from either of
$\Gamma$ or $\Gamma'$ by contracting the "small" triangle to the
point. This generates the mappings $i\colon H_1(\Gamma,
\mathbb{Z}_2) \to H_1(\Gamma'', \mathbb{Z}_2)$ and $i' \colon
H_1(\Gamma', \mathbb{Z}_2) \to H_1(\Gamma'', \mathbb{Z}_2)$, see
Fig. 4, right. Then the following holds: the cocycle $h$ is equal to
zero for small triangles, besides that if for $a\in  H_1(\Gamma,
\mathbb{Z}_2)$, $a'\in H_1(\Gamma', \mathbb{Z}_2)$ we have $i(a) =
i'(a')$, then $h(a) = h(a')$.

\end{enumerate}

\begin{figure}
\centering\includegraphics[width=150pt]{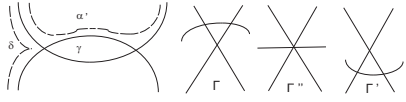} \caption{Second
and third moves}
\end{figure}

Let us consider now graphs with $n>1$ unicursal components.

{\bf Theorem 1.1}~In this case $H_1(\Gamma, \mathbb{Z}_2)$ is
generated not only by halves but by the following set of cycles:

1).~All halves,

2).~All bigons,

3).~(for $n\geqslant 3$) Cycles on $\Gamma$ corresponding to cycles
on the intersection graph of unicursal components.

\begin{proof}

Consider an arbitrary cycle $\gamma$ on $\Gamma$. Denote by
$M_{\Gamma}$ the set of cycles described in the statement of the
theorem and by Y the set of vertices where $\gamma$ rotates.

Let $\gamma$ rotates in a vertex $y\in Y$ which is incident to edges
$l_1, l_2\in \Gamma$ and the edges belong to the same unicursal
component. If $l_1$ and $l_2$ belong to the same half (say,
$\Gamma_{1,y}$), then by adding the half $\Gamma_{1,y}$ to the graph
$\Gamma$ (note than both halves are the elements of $M_{\Gamma}$) we
can obtain a cycle $\gamma'$ which rotates only in vertices from
$Y\setminus \{y\}$. If $l_1\in \Gamma_{1,y}$ and $l_2\in
\Gamma_{2,y}$ we can add any of the two halves to $\Gamma$ to
decrease the number of vertices in $Y$.

After applying this procedure several times obtained cycle
$\tilde{\gamma}$ will rotate only in those vertices where two pairs
of incident half-edges belong to other unicursal components. Thus
$\tilde{\gamma}$ belong to $M_{\Gamma}$.
\end{proof}

We use properties 2)-4) of $h$ to define it in the case of $n>1$.

{\bf Definition 1.6}~An element $h\in H^1(\Gamma, \mathbb{Z}_2)$ is
called {\it a homological parity} if the following conditions hold:

1).~$h(l_i)=0$ for each unicursal component $l_i$ of any framed
4-graph $\Gamma$.

2)-4). This properties coincide with ones listed above for $n=1$.

{\bf Definition 1.7}~We say that homological parity $h\in
H^1(\Gamma, \mathbb{Z}_2)$ {\it agrees} with a parity $\mathcal{P}$
if for any cycle $\gamma \in \Gamma$ holds $h(\gamma)= \sum
\limits_{v} \mathcal{P}(\gamma)$. Here the sum is over all vertices
$v$ such that $\gamma$ rotates in $v$.

{\bf Example 1.3}~Suppose that for each unicursal component $l_i$
$(i=1,\ldots,n)$ of a graph $\Gamma$ the total number of its
intersection points with other components is even. Denote the set of
all framed 4-graphs with this property by $\mathcal{G}$. We say that
homological parity is {\it Gaussian} if its value on any cycle
$\gamma$ is equal to the number of vertices where $\gamma$ go
transversally. Gaussian homological parity is well defined on the
set $\mathcal{G}$.

{\bf Definition 1.8}~A parity is called {\it Gaussian} if it agrees
with the Gaussian homological parity (in a sense of Definition 1.7).

Note that this definition of Gaussian parity coincides with one
given in Example 1.1.

\bigskip

The parity axioms listed at the beginning of this section lead to a
filtration on the set of 4-valent graphs and on the set of virtual
link diagrams.

Consider a virtual link diagram $L$. Let $f(L)$ be the diagram
obtained from $L$ by making all odd crossings virtual. By analogy
one can state the graph version of the action of $f$.

The next fact follows from definition.

{\bf Theorem 1.2}~The map $f$ is a well-defined map on the set of
all virtual links.

\begin{proof}
We must prove that if link diagrams $L_1$ and $L_2$ are equivalent
by a Reidemeister move then the diagrams $f(L_1)$ and $f(L_2)$ are
also equivalent. We check this only for the third Reidemeister move
(other cases are similar).

Recall that if vertices $a$, $b$, $c$ take part in the third move
than there is either one or three even vertices among them
(according to the third parity axiom). If all the three vertices are
even, then $f(L_1)=L_1$ and $f(L_2)=L_2$, so that $f(L_1)$ and
$f(L_2)$ also differ by the same move. In the other case virtual
diagrams are equivalent by {\it the detour move} (this is one of the
generalized Reidemeister moves).
\end{proof}

Note that analogous theorem holds for {\it any} knot-like objects
(e.g., for free links) and for {\it any} parity on these objects.

It turns out that one can construct a filtration on the set of all
virtual link by means of the map $f$. Let us introduce the following
notation: denote the set of link diagrams having all even crossings
by $\mathcal{A}^0$, the set of link diagrams $L$ such that $f(L)\in
\mathcal{A}^0$ by $\mathcal{A}^1$, etc. Than there is a natural
filtration $\mathcal{A}^0 \subset \mathcal{A}^1 \subset \ldots
\subset \mathcal{A}^n \subset \ldots$


\section{Atoms and their orientability}
\

The notion of atom was introduced by A.T.Fomenko in \cite{F} for the
study of bifurcations of integrable Hamiltonian systems. In
particular, atoms describe the structure of Morse functions on
2-manifolds in the neighbourhood of the critical level with several
critical points belonging to it. First connections between atoms and
knots were discovered by V.~O.~Manturov in papers on graphs embedded
into surfaces and his construction of Khovanov homology for virtual
links (see \cite{Ma2},\cite{Ma3}).

Let $P$ be a smooth closed compact 2-manifold. Let $\Gamma$ be a
graph of valency 4 embedded into P such that the graph splits $P$
into cells.

{\bf Definition 2.1}~A pair $(P, \Gamma)$ is called {\it an atom} if
the set of connected components in $P\setminus \Gamma$ can be
divided into two subsets (black and white cells) so that every edge
of $\Gamma$ is incident to one black and to one white cell.

{\bf Definition 2.2}~The graph $\Gamma$ is called {\it a frame} of
the atom $(P, \Gamma)$.

Atoms will be considered up to a natural isomorphism.

{\bf Definition 2.3}~Two atoms $(P_1, \Gamma_1)$ and $(P_2,
\Gamma_2)$ are called {\it isomorphic} if there exists a
homeomorphism $\phi\colon P_1 \to P_2$, taking frame into frame,
black cells to black cells, and white cells to white ones.

{\bf Definition 2.4}~An atom $(P, \Gamma)$ is called {\it
orientable} if the surface $P$ is orientable, otherwise it is called
non-orientable.

The framing of $\Gamma$ is induced by its embedding into $P$.

{\bf Definition 2.5}~One says that an orientation of edges of
$\Gamma$ determines {\it the source-sink structure} if they can be
oriented in such a way that in every vertex two opposite edges are
incoming, and the other two are emanating.

Note that if it is possible to determine the source-sink structure
on $\Gamma$ then it is determined completely by choosing orientation
of any single edge of $\Gamma$.

The following three theorems were proved in \cite{Ma4}:

{\bf Theorem 2.1}~An atom $(P, \Gamma)$ is orientable if and only if
it is possible to determine the source-sink structure on its frame.

{\bf Theorem 2.2}~An atom $(P, \Gamma)$ is orientable if and only if
any cycle of $\Gamma$ is orientable.

{\bf Theorem 2.3}~A cycle of an atom $(P, \Gamma)$ is orientable if
and only if this cycle goes through opposite half-edges in even
number of vertices of $\Gamma$.

We point out that starting from a frame with $n$ vertices it is
possible to construct $2^n$ atoms (some of them might turn out to be
isomorphic) choosing a pair of non-opposite half-edges in every
vertex, which form the border of black cells. It follows from
Theorem 2.1 that all atoms obtained in this way are either
simultaneously orientable or non-orientable.


\section{Proof of the main theorem}
\ {\bf Theorem 3.1}~If two virtual diagrams $K$ and $K'$ have
orientable atoms and they are equivalent by Reidemeister moves, then
there is a sequence of diagrams $K = K_1 \to \ldots \to K_n=K'$ all
having orientable atoms where $K_i$ is obtained from $K_{i-1}$ by a
Reidemeister move.

\begin{proof}
By standard applying of the second Reidemeister move one can obtain
another sequence $K_{1}\to \ldots \to K_{1}' \to K_{2}' \ldots \to
K_{n}'\to \ldots \to K_{n}$ where every two neighbouring diagrams
differ by a Reidemeister move and the intersection graph of
unicursal components is connected.

Since the diagram $K_{1}'$ has an orientable atom there exists the
source-sink structure on $K_{1}'$. Therefore $K_{1}'\in \mathcal{G}$
(here $\mathcal{G}$ is the set of framed 4-graphs defined in Example
1.3). For any unicursal component of $K_{1}'$ the parity of the
total number of intersection points with other components is an
invariant under Reidemeister moves. So we have
$K_{2}',\ldots,K_{n}'\in \mathcal{G}$.

Consider the Gaussian homological parity $h$ on $\mathcal{G}$.
Choose a Gaussian parity $\mathcal{P}_1$ such that
$\mathcal{P}_1(v)=0$ for any vertex $v\in K_{1}'$.

It follows from orientability of atoms for $K_{n}'$ that for each
vertex $w\in K_{n}'$ holds $\mathcal{P}_1(w)=0$.


Apply the map $f$ to the sequence of diagrams $K_{1}' \to \ldots \to
K_{n}'$ as many times as necessary to make all intermediate link
diagrams having orientable atom. Since $f(K_{1}')=K_{1}'$ and
$f(K_{n}')=K_{n}'$ we will obtain a sequence of diagrams with
orientable atoms between $K_{1}'$ and $K_{n}'$.

Then we perform the same method to sequences $K_1 \to \ldots \to
K_{1}'$ and $K_{n}' \to \ldots \to K_n$. In both cases we obtain
sequences with all intermediate diagrams having orientable atoms.

This completes the proof of the Theorem.
\end{proof}

\end{document}